\documentclass[reqno]{amsart}

\usepackage{amsmath,amssymb,amsthm,amsfonts,mathtools}
\usepackage{color,graphicx,tikz-cd,relsize,mathrsfs,latexsym}
\usepackage{hyperref,cleveref,enumerate,MnSymbol}

\makeatletter

\newtheorem{thm}{Theorem}[section]

\newtheorem{lem}[thm]{Lemma}

\newtheorem{prop}[thm]{Proposition}

\theoremstyle{definition}

\newtheorem{rem}[thm]{Remark}
\newtheorem{eg}[thm]{Example}

\numberwithin{equation}{section}

\allowdisplaybreaks

\title[Variations of cb maps on operator spaces]{Variations of completely bounded maps on operator spaces}
\author[Janson Antony]{Janson Antony}
\address{Department of Mathematics, University of Delhi, New Delhi - 110007, India}
\email{janson.math@gmail.com}
\author[Ajay Kumar]{Ajay Kumar$^\ast$}
\address{Department of Mathematics, University of Delhi, New Delhi - 110007, India}
\email{akumar@maths.du.ac.in}

\begin{document}

\begin{abstract} We introduce weighted cb maps and $\Lambda_\mu$-cb maps on operator spaces which are generalizations of completely bounded maps and a certain class of bilinear maps on operator spaces which we call $\lambda_\mu$-cb bilinear maps. Some basic properties of these maps and an operator space tensor product associated to $\lambda_\mu$-cb bilinear maps have been studied.
\end{abstract}
\thanks{$^\ast$Corresponding author, E-mail address: akumar@maths.du.ac.in}

\maketitle

\section*{Introduction}
\label{sec_intro-prelim}
An operator space, more precisely a concrete operator space is a (closed) subspace of a C$^*$-algebra. We refer to \cite{Ef_Ru_B} and \cite{Pis_B} for basic knowledge on operator spaces and operator space tensor products. 
Let $\mathcal{B}(\mathcal{H})$ denote the C$^*$-algebra of all bounded linear maps on a Hilbert space $\mathcal{H}$.
For two operator spaces $E,F$ which are subspaces of C$^*$-algebras $\mathcal{A}\subseteq \mathcal{B}(\mathcal{H}_1), \mathcal{B}\subseteq \mathcal{B}(\mathcal{H}_2)$ respectively, define the \emph{min} norm denoted by $\Vert\cdot\Vert_{min}$ on the algebraic tensor product $E\bigotimes F$ using the natural embeddings, 
$E\bigotimes F\subseteq \mathcal{A}\bigotimes \mathcal{B}\subseteq \mathcal{B}(\mathcal{H}_1)\bigotimes \mathcal{B}(\mathcal{H}_2)\subseteq \mathcal{B}(\mathcal{H}_1\bigotimes\mathcal{H}_2)$. 
Completion of the algebraic tensor product $E\bigotimes F$ under this norm, denoted by $E\bigotimes^{min} F$, is the minimal tensor product, which is often called as the operator space injective tensor product or the spatial tensor product. 
Clearly, $E\bigotimes^{min} F$ is an operator space as it is a closed subspace of the C$^*$-algebra $\mathcal{B}(\mathcal{H}_1\bigotimes\mathcal{H}_2)$.
For an operator space $E\subseteq \mathcal{B}(\mathcal{H})$, there are natural norms on $M_n(E)$, the space of all $n\times n$ matrices with entries from $E$, using the identification,
$M_n(E)=M_n\bigotimes^{min}E$
or equivalently using the embedding, 
$M_n(E)\subseteq M_n(\mathcal{B}(\mathcal{H}))=\mathcal{B}(\mathcal{H}^n)$, where $M_n$ is the C$^*$-algebra of all $n\times n$ matrices with scalar entries, identified with $\mathcal{B}(\mathbb{C}^n)$. 
This sequence of norms satisfies Ruan's axioms \cite{Ef_Ru_B}, and surprisingly, these norms on $E$ uniquely (unique in the sense described as in \cite{Ef_Ru_B}) determines the embedding $E\subseteq \mathcal{B}(\mathcal{H})$.

Ruan proved that any vector space $E$ together with a sequence of matrix norms $\Vert\cdot \Vert=(\Vert\cdot\Vert_n)_{n\in \mathbb{N}}$ where $\Vert\cdot\Vert_n$ is a matrix norm on $M_n(E)$ satisfying the \emph{Ruan's Axioms} \cite{Ru} can be embedded as a subspace of some C$^*$-algebra, say $\mathcal{A}_E$, such that the sequence of matrix norms induced by the inclusion $E\subseteq \mathcal{A}_E$ coincides with $\Vert\cdot\Vert$. An operator space defined in this way is commonly called an abstract operator space.

For operator spaces $E$ and $F$, a linear map $\phi:E\to F$ is said to be \emph{completely bounded} (cb in short) if the associated maps $\phi^{(n)}:M_n(E)\to M_n(F)$ defined as $\phi^{(n)}([e_{ij}]):=[\phi(e_{ij})]$ are uniformly bounded, where $n\in \mathbb{N}$, and is said to be a \emph{complete isometry} if the map $\phi_n$ turns out to be an isometry for every $n$.
Intuitively, one may think of cb maps as those which respect the matrix norms at each level.
The set of all cb maps denoted by $\mathcal{CB}(E,F)$ forms a normed linear space by defining the cb norm as $\Vert \phi\Vert_{cb}:=\sup_{n\in \mathbb{N}}\Vert \phi^{(n)}\Vert$, for $\phi\in \mathcal{CB}(E,F)$. 
It can further be given an operator space structure by identifying $M_n(\mathcal{CB}(E,F))$ with $\mathcal{CB}(E,M_n(F))$. This identification defines an operator space structure on the dual of an operator space $E$ when we choose $F=\mathbb{C}$.

There is no doubt that the most appropriate morphisms in the category of operator spaces are completely bounded maps. However, we can also have other special kind of linear maps between operator spaces, we shall discuss two of its kind in this article. 
In Section \ref{sec_weighted-cb}, we introduce weighted cb maps and under Section \ref{sec_lambda-cb}, $\Lambda_\mu$-cb maps. 
The sets of all weighted cb and $\Lambda_\mu$-cb maps, are given natural operator space structures.
In Section \ref{sec_bil}, for three operator spaces $E,F$, and $G$, we introduce a certain class of bilinear maps from $E\times F$ to $G$, which we call as completely $\lambda_\mu$-bounded bilinear maps, in a similar fashion as how completely bounded bilinear maps and jointly completely bounded bilinear maps are defined \cite{Ef_Ru,Bl_Pn}. 
With suitable choice of $\lambda$ and $\mu$, the set of all completely $\lambda_\mu$-bounded bilinear maps, denoted by $\mathrm{CB}_\lambda^\mu(E\times F,G)$ become an operator space in a natural way and we associate a tensor norm on the algebraic tensor product $E\bigotimes F$ in such a way that its dual become completely isometric to the space $\mathrm{CB}_\lambda^\mu(E\times F,\mathbb{C})$.

\section{Preliminaries and notations}
\label{sec_prelim-notations}
Let $E$ be a concrete operator space. Then the matrix norms on $M_n(E)$ satisfies the following two conditions:
\begin{enumerate}[\quad(R1)]
\item $\Vert e_1\oplus e_2\Vert_{n+m}\leq \max\{\Vert e_1\Vert_n,\Vert e_2\Vert_m\}$ for any $e_1\in M_n(E)$ and $e_2\in M_m(E)$. 
\item $\Vert \alpha e\beta\Vert_m\leq \Vert\alpha\Vert \Vert e\Vert_n \Vert \beta\Vert$ for any $e\in M_n(E)$, $\alpha\in M_{m\times n}$ and $\beta\in M_{n\times m}$.
\end{enumerate} 
We refer to the above two conditions as Ruans's Axioms.

Whenever $X,Y$ are normed linear space we shall denote by $\mathcal{L}(X,Y)$ and $\mathcal{B}(X,Y)$, respectively the spaces of all linear maps and bounded linear maps from $X$ to $Y$.
Let $E,F$ be any two operator spaces. 
Every completely bounded map from $E$ to $F$ is bounded, and the converse holds whenever $F$ is finite dimensional (or the map has finite rank) or if $F$ is a subspace of a commutative C$^*$-algebra. A well known example of a bounded map which is not cb is the usual involution on $\mathcal{B}(\ell^2)$ given by $T\mapsto T^*$ where $T^*$ denotes the adjoint of a bounded linear operator $T$ on $\ell^2$.

The smallest and the largest Banach space cross tensor norms on the algebraic tensor product of two Banach spaces are called respectively Banach space injective norm denoted by $\Vert\cdot\Vert_{\nu}$ and Banach space projective norm denoted by $\Vert\cdot\Vert_{\gamma}$ which are defined on the algebraic tensor product of two Banach spaces $X,Y$ as, $\Vert u\Vert_\nu=\sup \sum_{i=1}^n |f(x_i)g(y_i)|$ where $f\in X^*_1,g\in Y^*_1,u=\sum_{i=1}^n x_i\otimes y_i$ and $\Vert u\Vert_\gamma=\inf \{\sum_{i=1}^n\Vert x_i\Vert\Vert y_i\Vert\mid u=\sum_{i=1}^n x_i\otimes y_i\}$. 
We refer the reader to \cite{Ry_B,Bl} for all necessary details on tensor products of Banach spaces and operator algebras.
If $\Vert\cdot\Vert_\mu$ is a Banach space tensor norm, we shall always denote by $X\bigotimes_\mu Y$ the algebraic tensor product with this norm and its completion will be denoted by $X\bigotimes^\mu Y$.

Let $M_\infty(E)$ denote the set of all infinite matrices $[e_{ij}]_{1\leq i,j<\infty}$ with only finitely many non-zero entries from $E$. Clearly, $M_\infty(\mathbb{C})$ can naturally be embedded as a subspace of $\mathcal{K}(\ell^2)$, the C$^*$-algebra of all compact operators on $\ell^2$. Thus we may also identify $M_\infty(E)$ as $M_\infty(\mathbb{C})\bigotimes_{min} \mathcal{K}(\ell^2)$.
An element $u\in M_\infty(E\bigotimes F)$, where $E$ and $F$ are operator spaces, can be represented in three special ways as $u=\alpha(e\otimes f)\beta$, $u=\alpha(e\odot f)\beta$ and $u=\alpha(e\bullet f)\beta$ where $e=[e_{ij}]\in M_\infty(E),f=[f_{kl}]\in M_\infty(F),\alpha,\beta\in M_\infty(\mathbb{C})$ and $e\otimes f=[e_{ij}\otimes f_{kl}]$, $e\odot f=[\sum_{r=1}^\infty e_{ir}\otimes f_{rj}]$ and $e\bullet f=\alpha [e_{ij}\otimes f_{ij}]\beta$. 
One can define $\Vert u\Vert_\wedge=\inf\{\Vert\alpha\Vert\Vert e\Vert\Vert f\Vert\Vert\beta\Vert\mid u=\alpha(e\otimes f)\beta\}$, $\Vert u\Vert_h=\inf\{\Vert\alpha\Vert\Vert e\Vert\Vert f\Vert\Vert\beta\Vert\mid u=\alpha(e\odot f)\beta\}$ and $\Vert u\Vert_s=\inf\{\Vert\alpha\Vert\Vert e\Vert\Vert f\Vert\Vert\beta\Vert\mid u=\alpha(e\bullet f)\beta\}$, which gives three operator space structures to $E\bigotimes F$, whose completion with respect to these norms are named as the  operator space projective tensor product, the Haagerup tensor product and the Schur tensor product respectively \cite{AD_DW,VR_AK}. 
All necessary details on operator space tensor products can be seen in \cite{Bl,Ef_Ru_B} and \cite{Pis_B}.

An operator space tensor product $\bigotimes^\mu$ is said to be injective if whenever $\phi_1:E_1\to F_1$ and $\phi_2:E_2\to F_2$ are complete isometries, then, so is the tensor product map $\phi_1\otimes\phi_2:E_1\bigotimes_\mu E_2\to F_1\bigotimes_\mu F_2$. 
On the other hand, if $\phi_1\otimes \phi_2$ turn out to be completely bounded whenever $\phi_1$ and $\phi_2$ are completely bounded with $\Vert\phi_1\otimes\phi_2\Vert_{cb}\leq \Vert\phi_1\Vert_{cb}\Vert\phi_2\Vert_{cb}$, then we say $\mu$ is functorial. 
The operator space injective and the Haagerup tensor products are well known examples for injective tensor products while most of the tensor products that we consider including the operator space projective and the Schur tensor product are functorial \cite{Ef_Ru_B,VR_AK}.
An operator space tensor norm $\Vert\cdot\Vert_\mu$ is said to be matrix subcross if for any $e\in M_n(E)$ and $f\in M_n(F)$, $\Vert e\otimes f\Vert_\mu\leq \Vert e\Vert\Vert f\Vert$ and if equality holds then we call it matrix cross.

The elementary matrix in $\mathbb{M}_n$ whose $ij^\text{th}$ entry is 1 and all other entries are zeroes will be denoted by $\varepsilon_{ij}$ so that any $e=[e_{ij}]\in M_n(E)$ can be written as $e=\sum_{i,j}e_{ij}\otimes \varepsilon_{ij}$.

\section{Weighted cb maps on operator spaces}
\label{sec_weighted-cb}
Let $E,F$ be two operator spaces and $\lambda=(\lambda_n: M_n\to M_n)_{n\in \mathbb{N}}$ be a uniformly bounded sequence of non-zero bounded linear maps. 
A linear map $\phi:E\to F$ is said to be weighted completely bounded with $\lambda$ as the weight ($\lambda$-cb in short), if the associated maps $\phi\otimes \lambda_n: E\bigotimes_{min}M_n\to F\bigotimes_{min}M_n$ are uniformly bounded, i.e. $\sup_{n\in \mathbb{N}}\Vert\phi\otimes \lambda_n\Vert<\infty$. 
If $\phi$ is $\lambda$-cb, then we set $\Vert\phi\Vert_{cb}^\lambda=\sup_{n\in \mathbb{N}}\Vert\phi\otimes \lambda_n\Vert$.
The collection of all $\lambda$-cb maps from $E$ to $F$, denoted by $\mathcal{CB}_\lambda(E,F)$ is a linear subspace of $\mathcal{L}(E,F)$ and $\Vert\cdot\Vert_{cb}^\lambda$ is a norm on it. We leave the details to the reader.

A special kind of weights are those obtained by matrix conjugation by unitaries. 
If $U_n\in M_n$ are unitary matrices for each $n\in \mathbb{N}$, we may define $\lambda_n:M_n\to M_n$ as $\lambda_n(A)=U_n^{-1}AU_n$. 
In this case, if $\phi:E\to F$ is linear, then $\phi\otimes\lambda_n$ will be acting on $M_n(E)$ as $(\phi\otimes\lambda_n)(e)=U_n^{-1}\phi^{(n)}(e)U_n$ for any $e\in M_n(E)$. 
Then, $\mathcal{CB}_\lambda(E,F)=\mathcal{CB}(E,F)$ isometrically, because, for any $e\in M_n(E)$, if $U\in M_n$ is a unitary, then $\Vert U^{-1}eU\Vert=\Vert e\Vert$, which easily follows from Ruan's axiom R2.

Another interesting example is the choice of $\lambda_n$ as the transpose map on $M_n$ for every $n$. Here, $\lambda_n$ is an isometry for every $n$. 
Clearly, for any $\phi:E\to F$ and $e=[e_{ij}]\in M_n(E)$, we have $(\phi\otimes \lambda_n)(e)=[\phi(e_{ji})]$. 
Thus, it easily follows that the identity map on $\mathcal{K}(\ell^2)$ fails to be $\lambda$-cb, whereas the adjoint map succeeds to be $\lambda$-cb. 

\begin{prop} 
\label{prop_lambdacb_composition}
Let $E,F,G$ be operator spaces. If $\phi\in \mathcal{CB}(E,F)$ and $\psi\in \mathcal{CB}_\lambda(F,G)$, or $\phi\in \mathcal{CB}_\lambda(E,F)$ and $\psi\in \mathcal{CB}(F,G)$, then $\psi\circ\phi\in \mathcal{CB}_\lambda(E,G)$.
\begin{proof}
Let $e\in M_n(E)$. 
Consider
\[
((\psi\circ\phi)\otimes \lambda_n)(e)=\left((\psi\otimes \lambda_n)\circ \phi^{(n)}\right)(e)=\left(\psi^{(n)}\circ(\phi\otimes \lambda_n)\right)(e).
\]
Thus, if $\phi$ is cb and $\psi$ is $\lambda$-cb, then
\[
\Vert((\psi\circ\phi)\otimes \lambda_n)(e)\Vert=\left\Vert\left((\psi\otimes \lambda_n)\circ \phi^{(n)}\right)(e)\right\Vert\leq \Vert\psi\Vert_{cb}^\lambda\Vert\phi\Vert_{cb}\Vert e\Vert.
\]
Hence, $\psi\circ\phi$ is $\lambda$-cb with $\Vert\psi\circ\phi\Vert_{cb}^\lambda\leq \Vert\psi\Vert_{cb}^\lambda\Vert\phi\Vert_{cb}$.
A similar calculation works if we take $\phi$ as $\lambda$-cb and $\psi$ as cb.
\end{proof}
\end{prop}

As in the case of cb-maps, using the uniform boundedness of the sequence $\lambda$, we can impose conditions on the space $F$ so that any bounded linear map $\phi:E\to F$ become $\lambda$-cb for any $\lambda$. 

\begin{prop}
\label{prop_lambdacb}
Let $E,F$ be operator spaces and $\phi\in \mathcal{B}(E,F)$. Then, $\phi\in \mathcal{CB}_\lambda(E,F)$ if any of the following conditions is/are satisfied.
\begin{enumerate}[(i)]
\item\label{prop_lambdacb_functionals} $F=\mathbb{C}$, i.e. $\phi$ is a bounded linear functional.
\item\label{prop_lambdacb_finiterank} $F$ is finite dimensional or if $\phi$ is of finite rank.
\item\label{prop_lambdacb_commrange} $F$ is a subspace of some commutative C$^*$-algebra.
\end{enumerate}
\begin{proof}
\begin{enumerate}[(i)]
\item 
Consider $\phi\otimes \lambda_n:M_n(E)\to M_n(\mathbb{C})$ given by $e=[e_{ij}]\mapsto \sum_{i,j}\phi(e_{ij})\otimes \lambda_n(\varepsilon_{ij})=\lambda_n(\phi^{(n)}(e))$.
Hence
\[
\Vert (\phi\otimes \lambda_n)(e)\Vert=\Vert \lambda_n(\phi^{(n)}(e))\Vert\leq \sup_{n\in \mathbb{N}}\Vert \lambda_n\Vert\Vert\phi\Vert_{cb}\Vert e\Vert.
\]
Thus $\phi$ is $\lambda$-cb with $\Vert\phi\Vert_{cb}^\lambda\leq \sup_{n\in \mathbb{N}}\Vert \lambda_n\Vert\Vert\phi\Vert$.

\item 
Let $\{f_k\mid 1\leq k\leq n\}$ be an Auerbach basis for the range of $\phi$, i.e. $f_k$ are unit vectors and there exists bounded linear functionals $g_k$ for $1\leq k\leq n$ such that $g_i(f_j)=\delta_{i,j}$. Without loss of generality we can assume that $g_k$ is defined on whole of $F$.
For each $k$, the bounded linear functional $g_k\circ\phi$ being $\lambda$-cb, and the maps $\theta_k:\mathbb{C}\to F$ mapping $c\mapsto cf_k$ being cb, we can conclude from Proposition \ref{prop_lambdacb_composition} that $\theta_k\circ (g_k\circ \phi)$ is $\lambda$-cb and hence $\phi=\sum_{k=1}^n \theta_k\circ g_k\circ\phi$ is also $\lambda$-cb.

\item 
By the injectivity of the minimal tensor product, without loss of generality we can assume that $F=C(X)$, the C$^*$-algebra of all continuous scalar valued functions on a compact Hausdorff $X$. Moreover, we have the identification $C(X)\bigotimes_{min}M_n=C(X,M_n)$, where the latter is the C$^*$-algebra of all continuous $M_n$ valued functions on $X$.
Consider $\phi\otimes \lambda_n:M_n(E)\to C(X,M_n)$ given by $e=[e_{ij}]\mapsto f$ where $f:X\to M_n$ is the function $x\mapsto \lambda_n[\phi(e_{ij})(x)]=\lambda_n(\phi^{n}(e)(x))$. 
Hence 
\[
\Vert (\phi\otimes \lambda_n)(e)\Vert=\sup_{x\in X}\Vert \lambda_n(\phi^{n}(e)(x))\Vert\leq \Vert \lambda_n\Vert\Vert\phi\Vert_{cb}\Vert e\Vert\leq \sup_{n\in \mathbb{N}}\Vert \lambda_n\Vert\Vert\phi\Vert\Vert e\Vert.
\]
Thus $\phi$ is $\lambda$-cb with $\Vert\phi\Vert_{cb}^\lambda\leq \sup_{n\in \mathbb{N}}\Vert \lambda_n\Vert\Vert\phi\Vert$.\qedhere
\end{enumerate}
\end{proof}
\end{prop}

Now we shall give an operator space structure to $\mathcal{CB}_\lambda(E,F)$. 

\begin{lem} Let $E,F$ be operator spaces and $\phi_{ij}\in \mathcal{CB}_\lambda(E,F)$ for $1\leq i,j\leq n$. Then, the map $\phi:E\to M_n(F)$ defined as $\phi(e)=[\phi_{ij}(e)]$ is $\lambda$-cb.
\begin{proof}
Consider $\phi\otimes \lambda_n:E\bigotimes M_m\to M_n(F)\bigotimes M_m$.
For $e=[e_{kl}]=\sum_{k,l}e_{kl}\otimes \varepsilon_{kl}\in E\bigotimes M_m$, we have
\[
(\phi\otimes\lambda_m)(e)=\sum_{k,l}\phi(e_{kl})\otimes \lambda_m(\varepsilon_{kl})= \sum_{k,l}[\phi_{ij}(e_{kl})]_{i,j}\otimes \lambda_m(\varepsilon_{kl})=[(\phi_{ij}\otimes \lambda_m)(e)]_{ij}.
\]
Thus
\[
\Vert (\phi\otimes\lambda_m)(e)\Vert=\Vert [(\phi_{ij}\otimes \lambda_m)(e)]_{ij}\Vert\leq n^2\max_{i,j}\Vert\phi_{ij}\otimes\lambda_m\Vert\leq n^2\max_{i,j}\Vert\phi_{ij}\Vert_\lambda\Vert e\Vert,
\]
i.e. $\phi$ is $\lambda$-cb with $\Vert \phi\Vert_{cb}^\lambda\leq n^2\max_{1\leq i,j\leq n}\Vert\phi_{ij}\Vert_{cb}^{\lambda}$.
\end{proof}
\end{lem}

As in the case of $\mathcal{CB}(E,F)$, we can associate a sequence of matrix norms on $\mathcal{CB}_\lambda(E,F)$ using the identification $M_n(\mathcal{CB}_\lambda(E,F))=\mathcal{CB}_\lambda(E,M_n(F))$, and we expect that this sequence of matrix norms would give rise to an operator space structure on $\mathcal{CB}_\lambda(E,F)$. 

\begin{thm} $\mathcal{CB}_\lambda(E,F)$ is an operator space with matrix norms obtained from the identification $M_n(\mathcal{CB}_\lambda(E,F))=\mathcal{CB}_\lambda(E,M_n(F))$.
\begin{proof}
We shall verify the Ruan's Axioms.
Let $\phi_1\in M_n(\mathcal{CB}_\lambda(E,F))=\mathcal{CB}_\lambda(E,M_n(F))$ and $\phi_2\in M_m(\mathcal{CB}_\lambda(E,F))=\mathcal{CB}_\lambda(E,M_m(F))$. Then, for any $e\in E$, we have
\[
((\phi_1\oplus\phi_2)\otimes\lambda_k)(e)=((\phi_1\otimes\lambda_k)(e))\oplus((\phi_2\otimes\lambda_k)(e)),
\]
and $F$ being an operator space, we get
\begin{equation*}
\begin{aligned}
\Vert((\phi_1\oplus\phi_2)\otimes\lambda_k)(e)\Vert
&\leq\max\{\Vert(\phi_1\otimes\lambda_k)(e)\Vert,\Vert(\phi_2\otimes\lambda_k)(e)\Vert\}\\
&\leq\max\{\Vert(\phi_1\Vert_{cb}^\lambda \Vert e\Vert,\Vert\phi_2\Vert_{cb}^\lambda\Vert e\Vert\}.
\end{aligned}
\end{equation*}

Taking supremum over $e\in M_k(E)$ such that $\Vert e\Vert\leq 1$, we get required inequality as in R1.
Similarly, let $\phi=[\phi_{ij}]\in M_n(\mathcal{CB}_\lambda(E,F))$ and $\alpha,\beta\in M_n$. 
Let $e\in M_k(E)$. Consider
\[
(\alpha\phi\beta\otimes \lambda_k)(e)=\alpha((\phi\otimes \lambda_k)(e))\beta,
\]
and $F$ being an operator space, we get
\begin{equation*}
\begin{aligned}
\Vert(\alpha\phi\beta\otimes \lambda_k)(e)\Vert 
&\leq \Vert\alpha\Vert \Vert(\phi\otimes \lambda_k)(e)\Vert\Vert\beta\Vert\\
&\leq \Vert\alpha\Vert \Vert\phi\Vert_{cb}^\lambda\Vert\beta\Vert \Vert e\Vert.
\end{aligned}
\end{equation*}
Thus we have the required inequality as in R2.
\end{proof}
\end{thm}

Notice that, from Proposition \ref{prop_lambdacb}(\ref{prop_lambdacb_functionals}), we can conclude that the bounded linear functionals on an operator space coincides with the $\lambda$-cb functionals. We may even choose all $\lambda_n$ as isometries so that we have isometrically $\mathcal{CB}_\lambda(E,\mathbb{C})=\mathcal{CB}(E,\mathbb{C})$ for any operator space $E$.
But we still can not conclude that $\mathcal{CB}_\lambda(E,M_n(\mathbb{C}))=\mathcal{CB}(E,M_n(\mathbb{C}))$. 
As a result, the operator space structure on $\mathcal{CB}_\lambda(E,\mathbb{C})$ need not be the same as that on $\mathcal{CB}(E,\mathbb{C})$. 
Hence we shall denote the dual ($\lambda$-dual) space $\mathcal{CB}_\lambda(E,\mathbb{C})$ as $E^*_\lambda$ to distinguish it from the usual operator space dual.

\subsection*{Quantizations through tensor products}
\label{subsec_quant}
The process of defining an operator space structure on a given normed linear space $X$, either as an isometric embedding into a C$^*$-algebra or by explicitly defining a sequence of matrix norms satisfying the Ruan's axioms with the condition that the matrix norm on $M_1(X)$ coincides with the given norm on $X$, is called \emph{quantization}. 
It is to be noted that any normed linear space can be embedded naturally inside a commutative C$^*$-algebra, namely the C$^*$-algebra of all continuous complex valued functions on the closed unit ball of the dual space $X^*$ which is compact with respect to the weak$^*$ topology. This is popularly known as the MIN quantization. It is minimal in the sense that the matrix norms obtained by any other quantization on $X$ would be bigger than that obtained by the MIN quantization. 
On the other hand, there is a MAX quantization which is maximal in the sense analogues to the one mentioned above.
In fact, we can explicitly construct the matrix norms in the two cases as described in \cite[Section 3.3]{Ef_Ru_B}. We shall denote by $X_{MIN}$ and $X_{MAX}$ respectively the space $X$ with the MIN and MAX operator space structures, and the operator space norms on them shall be denoted by $\Vert\cdot\Vert_{MIN}$ and $\Vert\cdot\Vert_{MAX}$ respectively.

For a normed linear space $X$, the MIN quantization can be obtained as follows: For any $x\in M_n(X)=M_n\bigotimes X$, $\Vert x\Vert_{MIN}=\Vert x\Vert_{M_n\bigotimes^\nu X}$ where $\bigotimes^\nu$ denotes the Banach space injective tensor product.
Motivated by this, we ask the following question: If $X$ is a normed linear space and $\bigotimes^\mu$ denotes a Banach space tensor product, then can we define an operator space norm on $X$ by identifying $M_n(X)$ with $M_n\bigotimes_\mu X$? 
We shall only concentrate on tensor products which give rise to tensor norms lying in between the Banach space injective and the Banach space projective tensor norms, which would imply that the tensor norm is a cross norm. 
Of course, we shall consider the case when $X$ is an operator space and $\bigotimes^\mu$ is an operator space tensor product so that we may obtain various quantizations of the underlying normed linear space. 
In any case, for any $X$, the tensor product $M_n\bigotimes^\mu X$ being bicontinuously isomorphic to $M_n\bigotimes^\lambda X$, we can have the identification $M_n\bigotimes^\mu X=M_n(X)$ as vector spaces. 
Hence, we shall denote by $M_n^\mu(X)$ the space $M_n\bigotimes^\mu X$.
For notational convenience we shall use the phrase `$\bigotimes^\mu$ quantizes' or `$\bigotimes^\mu$ is a quantizing tensor product' to mean that for every normed linear space $X$, the identification $M_n(X)=M_n^\mu(X)$ defines an operator space structure on $X$. If $X$ is an operator space and $\bigotimes^\mu$ is an operator space tensor product, then we hope that it must be clear from the context that we are talking about a new (not necessarily different) operator space structure on $X$.

When we consider an injective operator space tensor product $\bigotimes^\mu$, the operator space $M_n$ being embedded completely isometrically inside $\mathcal{K}(\ell^2)$, we have $M_n\bigotimes^\mu E\subset \mathcal{K}(\ell^2)\bigotimes^\mu E$ completely isometrically for any operator space $E$ and for any $n\in \mathbb{N}$. 
Therefore, for any $e=[e_{ij}]\in M_n^\mu(E)$, we have $\Vert e\Vert_{M_n^\mu(E)}=\Vert e\Vert_{\mathcal{K}(\ell^2)\bigotimes^\mu E}$. 
Similarly, for $n,m\in \mathbb{N}$, we have completely isometric embeddings $M_n^\mu(E)\hookrightarrow M_{n+m}^\mu(E)$ and $M_m^\mu(E)\hookrightarrow M_{n+m}^\mu(E)$ given by 
\[e_1\mapsto \begin{bmatrix}
e_1 & 0\\
0 & 0
\end{bmatrix},
\quad
e_2\mapsto \begin{bmatrix}
0 & 0\\
0 & e_2
\end{bmatrix},
\]
for any $e_1\in M_n^\mu(E)$ and $e_2\in M_m^\mu(E)$.
Thus we have the completely isometric embedding $M_n^\mu(E)\oplus M_m^\mu(E)\subseteq M_{n+m}^\mu(E)$, which gives $\Vert e_1\oplus e_2\Vert=\max\{\Vert e_1\Vert,\Vert e_2\Vert\}$. This proves that matrix norms satisfies the Ruan's axiom R1 whenever $\mu$ is injective. 
We could not find a sufficient condition on $\bigotimes^\mu$ such that R2 may also be satisfied. 

Using a Banach space tensor product $\bigotimes^\mu$, one may extend the idea of weighted cb maps on operator spaces to those between normed linear spaces as follows.
Let $X,Y$ be two normed linear spaces and $\lambda=(\lambda_n)_{n\in \mathbb{N}}$ be a uniformly bounded sequence of linear maps where $\lambda_n:M_n\to M_n$. 
A linear map $\phi:X\to Y$ is said to be $\lambda_\mu$ completely bounded or $\lambda_\mu$-cb in short, if the associated maps $\phi\otimes \lambda_n:X\bigotimes^\mu M_n\to Y\bigotimes_\mu M_n$ are uniformly bounded.
But, if $\bigotimes^\mu$ is a quantizing tensor product, one may observe that it is not much different from the weighted cb maps, for, if $X', Y'$ denote the operator spaces obtained by quantizing $X$ and $Y$ respectively using $\bigotimes^\mu$, then it is easy to observe that the notion of $\lambda_\mu$-cb maps from $X\to Y$ coincides with the definition $\lambda$-cb maps from $X'$ to $Y'$.
Thus, in this case, it suffices to study properties of weighted cb maps on operator spaces to understand $\lambda_\mu$-cb maps on normed linear spaces. As a particular case, if we take the Banach space injective tensor product $\bigotimes^\nu$ in place of $\bigotimes^\mu$, then the operator spaces $X'$ and $Y'$ are respectively $X_{MIN}$ and $Y_{MIN}$. Hence, from Proposition \ref{prop_lambdacb}(\ref{prop_lambdacb_commrange}), $Y_{MIN}$ being subspace of a commutative C$^*$-algebra, we have $\mathcal{CB}_\lambda(X_{MIN},Y_{MIN})=\mathcal{B}(X_{MIN},Y_{MIN})=\mathcal{B}(X,Y)$.

If $\bigotimes^\mu$ is not a quantizing tensor product, then one may consider the notion of $\lambda_\mu$-cb maps defined above on normed linear spaces or on operator spaces which may give strange results. We do not plan to discuss it here. We shall consider a bilinear analogue of $\lambda_\mu$-cb maps on operator spaces under Section \ref{sec_bil}.

\section{$\Lambda_\mu$-cb maps on operator spaces}
\label{sec_lambda-cb}
Let $\Lambda$ be a collection of operator spaces containing atleast one non-trivial element (always assumed from here onwards). The collection may not be a set, we will be considering examples of $\Lambda$ being proper classes, for example, the class of all commutative C$^*$-algebras. 
Let $\bigotimes^\mu$ be an operator space tensor product not necessarily quantizing, but we shall assume that $\bigotimes^\mu$ is matrix cross, meaning that, for any two operator spaces $E,F$, we must have $\Vert e\otimes f\Vert_\mu=\Vert e\Vert\Vert f\Vert$ for all $e\in M_n(E)$ and $f\in M_m(F)$.
Let $E,F$ be operator spaces. A bounded linear map $\phi:E\to F$ is said to be \emph{$\Lambda_\mu$-completely bounded} ($\Lambda_\mu$-cb in short) if the associated maps $\phi_X^\mu:=\phi\otimes_\mu I_X:E\bigotimes_\mu X\to F\bigotimes_\mu X$ are uniformly bounded (i.e. $\sup_{X\in \Lambda}\Vert \phi_X^\mu\Vert<\infty$.), where $I_X$ denotes the identity map from $X$ to $X$. 
Denote by $\mathcal{CB}_\Lambda^\mu(E,F)$, the set of all $\Lambda_\mu$-cb maps from $E$ to $F$. 
When $E=F$, $\mathcal{CB}_\Lambda^\mu(E,E)$ may be denoted as $\mathcal{CB}_\Lambda^\mu(E)$. 
Define 
\[
\Vert\cdot\Vert_{cb}^{\Lambda_\mu}:\mathcal{CB}_\Lambda^\mu(E,F)\to \mathbb{R}\text{\quad as\quad }\Vert \phi\Vert_{cb}^{\Lambda_\mu}:=\sup_{X\in\Lambda}\Vert \phi_X^\mu\Vert.
\]
A routine verification shows the following.

\begin{prop}
For any two operator spaces $E,F$ and any collection $\Lambda$ of operator spaces, $\mathcal{CB}_\Lambda^\mu(E,F)$ forms a vector subspace of $\mathcal{L}(E,F)$ and the function $\Vert\cdot\Vert_{cb}^{\Lambda_\mu}$ defines a norm on it.
\end{prop}

In the very special case when $\mu$ is taken as the \emph{min} tensor product, we shall omit the symbol $\mu$ from our notations and call a $\Lambda_\mu$-cb map simply as a $\Lambda$-cb map, $\mathcal{CB}_\Lambda^\mu(E,F)$ as $\mathcal{CB}_\Lambda(E,F)$, $\phi_X^\mu$ as $\phi_X$ etc.
When the tensor product $\mu$ under consideration is injective, as in the cases of the \emph{min} and Haagerup tensor products, we may take $\Lambda$ as a collection of unital C$^*$-algebras, because, if $\Lambda=\{X_\lambda\}$ is a collection of operator spaces, then $\Vert \sum_{i=1}^n e_i\otimes x_i\Vert_{E\bigotimes_\mu X_\lambda}=\Vert \sum_{i=1}^n e_i\otimes x_i\Vert_{E\bigotimes_\mu \mathcal{A}_\lambda}$ for any element $\sum_{i=1}^n e_i\otimes x_i\in E\bigotimes X_\lambda$ and for every $\lambda\in \Lambda$, where $\mathcal{A}_\lambda$ is a unital C$^*$-algebra such that $X_\lambda\subseteq \mathcal{A}_\lambda$ completely isometrically.

Injectivity of the tensor product $\bigotimes^\mu$ simplify things. However, very important tensor products such as the operator space projective tensor product and the Schur tensor product are not injective in general, though, it may be functorial which is the second best thing that we can ask for, after injectivity. 

\begin{prop}
\label{prop_bdd}
Let $E,F$ be two operator spaces and $\Lambda$ be a collection of operator spaces. Then the following statements hold.
\begin{enumerate}[(i)]
\item For any $\phi\in \mathcal{CB}_\Lambda^\mu(E,F)$, we have, $\Vert \phi\Vert\leq \Vert \phi\Vert_{cb}^{\Lambda_\mu}$.
\item If $\bigotimes^\mu$ is functorial, then for any $\phi\in \mathcal{CB}(E,F)$, we have, $\phi\in \mathcal{CB}_\Lambda^\mu(E,F)$ and $\Vert \phi\Vert_{cb}^{\Lambda_\mu}\leq\Vert \phi\Vert_{cb}$.
\end{enumerate} 
\begin{proof}
\begin{enumerate}[(i)]
\item 
Let $X\in \Lambda$ be a non-trivial operator space and $x_0\in X$ be a unit vector. 
Consider $\phi_X^\mu:E\bigotimes^\mu X\to F\bigotimes^\mu X$ defined on rank one tensors as $\phi_X^\mu(e\otimes x)=\phi(e)\otimes x$, for all $e\in E$ and $x\in X$. 
Let $(e_n)$ be a sequence of unit vectors in $E$ such that $\Vert\phi(e_n)\Vert$ converges to $\Vert\phi\Vert$. 
Consider the sequence $\phi_X^\mu(e_n\otimes x_0)$. 
As $\Vert e_n\otimes x_0\Vert=1$ and $\Vert\phi_X^\mu(e_n\otimes x_0)\Vert=\Vert\phi(e_n)\otimes x_0\Vert=\Vert\phi(e_n)\Vert$, we can conclude that $\Vert\phi\Vert\leq\Vert\phi_X^\mu\Vert$ and hence $\Vert\phi\Vert\leq\Vert\phi\Vert_{cb}^{\Lambda_\mu}$.
\item 
The tensor product $\bigotimes^\mu$ being functorial, $\phi:E\to F$ and $I_X:X\to X$ being completely bounded, $\phi_X^\mu$ is also completely bounded with $\Vert\phi_X^\mu\Vert\leq \Vert\phi\Vert_{cb}\Vert I_X\Vert_{cb}=\Vert\phi\Vert_{cb}$, i.e. $\Vert\phi\Vert_{cb}^{\Lambda_\mu}=\sup_{X\in \Lambda}\Vert\phi_X^\mu\Vert\leq \Vert\phi\Vert_{cb}$. \qedhere
\end{enumerate}
\end{proof}
\end{prop}
From the above proposition, we have, $\mathcal{CB}(E,F)\subseteq \mathcal{CB}_\Lambda^\mu(E,F)$ whenever $\mu$ is injective (the set inclusion must not be confused with isometric embedding of normed linear spaces), in particular we have $\mathcal{CB}(E,F)\subseteq\mathcal{CB}_\Lambda(E,F)$ by taking $\mu=min$. 
But, of course, the containment can be strict, depending on the choice of $\Lambda$. 
For example, let $\Lambda\subset \{M_n\mid n\in \mathbb{N}\}$ be finite, and choose $E,F,\phi$ such that $\phi$ is bounded but not completely bounded. An easy way to do this is by considering the identity map $I:E_{MIN}\to E_{MAX}$ where $E$ is an infinite dimensional space and $E_{MIN}$ and $E_{MAX}$ denote respectively the minimal and the maximal quantizations of $E$.
On the other hand, if we choose $\Lambda=\{M_n\mid n\in \mathbb{N}\}$, then the definition of $\Lambda$-cb maps coincides with that of cb maps. Thus we have $\mathcal{CB}(E,F)=\mathcal{CB}_\Lambda(E,F)$ isometrically, in this case. 

As a simple consequence of the injectivity of the minimal tensor product, if $X,Y$ are operator spaces with $i_Y:Y\to X$ being a complete isometry, then for any operator space $E$ with $I:E\to E$ being the identity map, the tensor map $I\otimes i_Y:E\bigotimes_{min}Y\to E\bigotimes_{min}X$ is also a complete isometry. 
In other words, if $Y\subseteq X$, then for any $u=\sum_{i=1}^n e_i\otimes y_i$ where $e_i\in E$ and $y_i\in Y$, we have, $\Vert u\Vert_{E\bigotimes_{min}Y}=\Vert u\Vert_{E\bigotimes_{min}X}$.
Thus, if $\Lambda$ has the property that, for every $n\in \mathbb{N}$, there is a C$^*$-algebra $\mathcal{A}\in\Lambda$ and a completely isometric embedding of $M_n$ into $\mathcal{A}$, then every $\Lambda$-cb map is cb. 
This is because, if we choose $\mathcal{A}_n$ to be a C$^*$-algebra in $\Lambda$ satisfying $M_n\subseteq \mathcal{A}_n$, then by injectivity of the minimal tensor product, if $\phi:E\to F$ is a bounded map, then $\Vert\phi_{\mathcal{A}_n}\Vert\geq \Vert\phi^{(n)}\Vert$ where $\phi^{(n)}$ denotes the $n^\text{th}$ amplification of $\phi$. Thus taking supremum over $\Lambda$, we get, $\Vert\phi\Vert_{cb}^\Lambda\geq \Vert\phi\Vert_{cb}$.
Choosing $\Lambda=\{\mathcal{K}(\ell^2)\}$ gives an example for the above case, as $\mathbb{M}_n$ is a $*$-subalgebra of $\mathcal{K}(\ell^2)$ for every $n$. Another similar example is taking $\Lambda=\{M_n(\mathcal{A})\mid n\in \mathbb{N}\}$ where a C$^*$-algebra $\mathcal{A}$ is fixed. 

As every cb map is $\Lambda$-cb, we have aplenty of examples of $\Lambda$-cb maps no matter whatever $\Lambda$ is.  
If $F$ is a finite dimensional operator space or a commutative C$^*$-algebra, then any bounded map $\phi:E\to F$ is $\Lambda$-cb for any operator space $E$. 
In particular, any bounded linear functional $\phi$ on $E$ is $\Lambda$-cb with $\Vert\phi\Vert_{cb}^\Lambda=\Vert\phi\Vert$.
But in other cases, such as $F$ being finite dimensional, eventhough all cb maps are $\Lambda$-cb, their norms can be different, meaning that the identity map between the spaces $\mathcal{CB}(E,F)$ and $\mathcal{CB}_\Lambda(E,F)$ may not be an isometry. This will be discussed later when we consider two different operator space structures on the dual space.
From Proposition \ref{prop_bdd}, we conclude that, for $\phi\in \mathcal{CB}(E,F)$, we have, $\Vert\phi\Vert\leq\Vert\phi\Vert_{cb}^\Lambda\leq\Vert\phi\Vert_{cb}$ for any non-empty $\Lambda$ and the examples discussed above with slight modifications give strict inequalities (See Example \ref{eg_ckmn} below). Another immediate observation is that if $\Lambda_1\subset\Lambda_2$, then $\mathcal{CB}_{\Lambda_2}(E,F)\subset \mathcal{CB}_{\Lambda_1}(E,F)$, and for $\phi\in \mathcal{CB}_{\Lambda_2}(E,F)$, we have, $\Vert\phi\Vert_{cb}^{\Lambda_1}\leq \Vert\phi\Vert_{cb}^{\Lambda_2}$.

It is a clear fact that for $n>1$, there is no commutative C$^*$-algebra $\mathcal{A}$ such that $M_n$ is embedded as a subalgebra of $\mathcal{A}$. Thus it is quite natural to expect that there is no completely isometric embedding of $M_n$ into $\mathcal{A}$ too. Curiously, analogous conclusion about non-commutative C$^*$-algebras may fail. There exists non-commutative C$^*$-algebras which do not contain any $M_n$ as subalgebras, though which contain $M_n$ embedded as subspaces for some $n>1$, completely isometrically. The reduced C$^*$-algebra of the free group on two generators is such an example.

\begin{prop}
\label{commutativecsa}
Let $\Lambda$ be the class of all commutative (unital) C$^*$-algebras and $E,F$ be any two operator spaces. Then $\mathcal{CB}_\Lambda(E,F)=B(E,F)$ isometrically.
\begin{proof}
Let $\mathcal{A},\mathcal{B}$ be C$^*$-algebras such that $E,F$ are embedded into $\mathcal{A},\mathcal{B}$ completely isometrically, respectively. 
Let $\Lambda=\{C(K_i)\}_{i\in I}$, where $C(K_i)$ denotes the C$^*$-algebra of all continuous complex valued functions on a compact Hausdorff space $K_i$ with the sup norm denoted by $\Vert\cdot\Vert_\infty$ with pointwise algebraic operations. 
We have the identification $\mathcal{A}\bigotimes^{min}C(K_i)=C(K_i,\mathcal{A})$ where $C(K_i,\mathcal{A})$ denotes the C$^*$-algebra of all $\mathcal{A}$-valued continuous functions on $K_i$. 
Hence, we have the completely isometric embeddings, $E\bigotimes^{min}C(K_i)\subset C(K_i,\mathcal{A})$ and $F\bigotimes^{min}C(K_i)\subset C(K_i,\mathcal{B})$. 
Let $\phi:E\to F$ be a bounded linear map. Consider $\phi_{\mathcal{A}_i}:C(K_i,\mathcal{A})\to C(K_i,\mathcal{B})$. 
For $f\in F\bigotimes^{min}C(K_i)$ identified as a continuous $\mathcal{A}$ valued function on $K_i$, its image $\phi_i(f)$ is the continuous $\mathcal{B}$ valued function on $K_i$ defined as $\phi_i(f)(x_i)=\phi(f(x_i))$ for every $x_i\in K_i$.
It follows that $\Vert \phi_i\Vert=\sup_{\Vert f\Vert_\infty\leq 1,x_i\in X_i}\Vert \phi(f(x_i))\Vert=\Vert \phi\Vert$.
In fact, using the well known identity $E\bigotimes^{min} C(K)=E\bigotimes^\nu C(K)$ (isometrically), we could have easily arrived at the conclusion, without any further calculations, that $\Vert \phi_i\Vert=\Vert \phi\Vert\;\forall\; i$, where $\nu$ denotes the Banach space injective tensor product.
\end{proof}
\end{prop}

\begin{eg}
\label{eg_ckmn}
Fix $n\in \mathbb{N}$. Let $\Lambda_n$ be the collection of all closed $*$-subalgebras of $C(K,\mathbb{M}_n)$ for some compact Hausdorff space $K$. 
Let $\phi:E\to F$ be a bounded linear map between operator spaces. 
As $C(K,\mathbb{M}_n)$ can be identified with $C(K)\bigotimes^{min}M_n$, we have naturally the completely isometric identification, $E\bigotimes^{min}C(K,\mathbb{M}_n)=M_n(E\bigotimes^{min}C(K))=M_n(E)\bigotimes^{min} C(K)$. 
Hence the map $\phi_{C(K,\mathbb{M}_n)}$ is nothing but the map $(\phi^{(n)})_{C(K)}$ where $\phi^{(n)}$ is the $\text{n}^{\text{th}}$ amplification of $\phi$. 
But, from Proposition \ref{commutativecsa}, we have, $\Vert \phi\Vert=\Vert\phi_{C(K)}\Vert$ for every $\phi$, and hence, $\Vert\phi_{C(K,\mathbb{M}_n)}\Vert=\Vert\phi^{(n)}\Vert\leq n\Vert\phi\Vert$. 
Thus, $\mathcal{CB}_{\Lambda_n}(E,F)$ is bicontinuously isomorphic to $B(E,F)$ via the identity map.
\end{eg}

\begin{prop} The normed linear space $\mathcal{CB}_\Lambda^\mu(E,F)$ is complete whenever $F$ is complete.
\begin{proof}
Let $(\phi^n)_{n\in \mathbb{N}}$ be a Cauchy sequence in $\mathcal{CB}_\Lambda^\mu(E,F)$. 
As the norm on $B(E,F)$ is smaller than that on $\mathcal{CB}_\Lambda^\mu(E,F)$ for any $\Lambda$-cb map, we have that $(\phi^n)_{n\in \mathbb{N}}$ is Cauchy in $B(E,F)$ also. 
But, we know that $B(E,F)$ is complete whenever $F$ is complete. 
Hence $(\phi^n)_{n\in \mathbb{N}}$ converges to some $\phi\in B(E,F)$. 
We shall prove that $(\phi^n)_{n\in \mathbb{N}}$ converges to $\phi$ in $\Vert\cdot\Vert_{cb}^{\Lambda_\mu}$ norm and that $\phi\in \mathcal{CB}_\Lambda^\mu(E,F)$. 
Let $\varepsilon>0$ be given. 
As $(\phi^n)_{n\in \mathbb{N}}$ is Cauchy, there exists $N\in \mathbb{N}$ and $K>0$ such that 
$\Vert\phi^m-\phi^n\Vert_{cb}^{\Lambda_\mu}<\frac{\varepsilon}{2}$ for all $m,n\geq N$ and 
$\Vert\phi^n\Vert_{cb}^{\Lambda_\mu}<K$ for every $n\in \mathbb{N}$. 
As $\Vert\phi-\phi^n\Vert$ converges to zero, $\Vert\phi_X-\phi_X^n\Vert$ also converges to zero for every $X\in\Lambda$ as $n$ tends to infinity. 
Hence there exists $N_X\in \mathbb{N}$ such that $\Vert\phi_X-\phi_X^m\Vert<\frac{\varepsilon}{2}$ for every $m\geq N_X$. 
Let $X\in\Lambda$ and $m>\max\{N_X,N\}$. 
Consider $\Vert\phi_X-\phi_X^n\Vert\leq \Vert\phi_X-\phi_X^m\Vert+\Vert\phi_X^m-\phi_X^n\Vert<\frac{\varepsilon}{2}+\Vert\phi^m-\phi^n\Vert_{cb}^{\Lambda_\mu}<\varepsilon$ for every $n\geq N$. 
Hence $\Vert\phi-\phi^n\Vert_{cb}^{\Lambda_\mu}\leq \varepsilon$ for every $n\geq N$,
i.e. $(\phi^n)_{n\in \mathbb{N}}$ converges to $\phi$ in $\Vert\cdot\Vert_{cb}^{\Lambda_\mu}$ norm (inside $B(E,F)$).
Let $n\geq N$. 
Consider $\Vert \phi_X\Vert\leq \Vert\phi_X-\phi_X^n\Vert+\Vert\phi_X^n\Vert\leq \Vert(\phi-\phi^n)_X\Vert+\Vert\phi_X^n\Vert\leq \Vert\phi-\phi^n\Vert_{cb}^{\Lambda_\mu}+\Vert\phi^n\Vert_{cb}^{\Lambda_\mu}<\varepsilon+K$. 
Hence $\Vert\phi\Vert_{cb}^{\Lambda_\mu}\leq \varepsilon+K$ and thus $\phi\in \mathcal{CB}_\Lambda^\mu(E,F)$.
\end{proof}
\end{prop}

\begin{prop}
\label{banach-algebra}
Let $E,F,G$ be operator spaces. If $\phi\in \mathcal{CB}_\Lambda^\mu(E,F)$ and $\psi\in \mathcal{CB}_\Lambda^\mu(F,G)$, then $\psi\circ\phi\in \mathcal{CB}_\Lambda^\mu(E,G)$. Moreover, $\Vert\psi\circ\phi\Vert_{cb}^{\Lambda_\mu}\leq\Vert\psi\Vert_{cb}^{\Lambda_\mu}\Vert\phi\Vert_{cb}^{\Lambda_\mu}$. In particular, $\mathcal{CB}_\Lambda^\mu(E)$ is a Banach algebra.
\begin{proof}
Let $X\in\Lambda$. 
Consider $(\psi\circ\phi)_X^\mu:E\bigotimes_\mu X\to G\bigotimes_\mu X$ defined as $(\psi\circ\phi)_X(e\otimes x)=(\psi\circ\phi)(e)\otimes x$ on rank one tensors. 
Thus, $(\psi\circ\phi)_X=\psi_X\circ\phi_X$ and hence 
$\Vert (\psi\circ\phi)_X\Vert\leq \Vert\psi_X\Vert \Vert\phi_X\Vert\leq \Vert\psi\Vert_{cb}^{\Lambda_\mu} \Vert\phi\Vert_{cb}^{\Lambda_\mu}$.
As the above inequality is true for all $X\in\Lambda$, we can easily conclude that $\Vert\psi\circ\phi\Vert_{cb}^{\Lambda_\mu}\leq\Vert\psi\Vert_{cb}^{\Lambda_\mu}\Vert\phi\Vert_{cb}^{\Lambda_\mu}$.
\end{proof}
\end{prop}

In the special case as in Example \ref{eg_ckmn}, we can talk about an involution for $\mathcal{CB}_{\Lambda_n}(E,F)$. Because for any $\phi\in \mathcal{CB}_{\Lambda_n}(E,F)$, the usual adjoint operator $\phi^*:F^*\to E^*$ turns out to be a $\Lambda_n$-cb map with $\Vert\phi^*\Vert_{cb}^{\Lambda_n}=\Vert\phi^{*(n)}\Vert=\Vert\phi^{(n)}\Vert=\Vert\phi\Vert_{cb}^{\Lambda_n}$, i.e. $\mathcal{CB}_{\Lambda_n}(E,F)$ has an isometric involution. The involution in general case is not clear.

\begin{rem}
In Proposition \ref{banach-algebra}, if we relax the condition of both $\phi$ and $\psi$ being $\Lambda_\mu$-cb to only one of the maps being $\Lambda_\mu$-cb and the other being bounded, then clearly there is no guarantee that their composition is $\Lambda_\mu$-cb. 
Because, any bounded map $\phi$ on an operator space is a composition of a $\Lambda_\mu$-cb map and a bounded map, namely the identity map and $\phi$ itself! 
We have been working with a fixed $\Lambda$ and a fixed $\mu$. 
But, of course, we can find particular instances such as some triples $((\Lambda_1,\mu_1),(\Lambda_2,\mu_2),(\Lambda_3,\mu_3))$ such that $\psi\circ\phi\in \mathcal{CB}_{\Lambda_3}^{\mu_3}(E,G)$ whenever $\phi\in \mathcal{CB}_{\Lambda_1}^{\mu_1}(E,F)$ and $\psi\in \mathcal{CB}_{\Lambda_2}^{\mu_2}(F,G)$. We do not plan to discuss it here.  
\end{rem}

\subsection{An operator space structure on $\mathcal{CB}_\Lambda^\mu(E,F)$}
\label{subsec_lambda-cb_oss}
Now we shall associate a sequence of matrix norms with $\mathcal{CB}_\Lambda^\mu(E,F)$, satisfying the Ruan's Axioms R1 and R2, so that it become an operator space. 
We will make use of the following lemma.

\begin{lem} If $E,F$ are operator spaces and $\phi_{ij}:E\to F$ are $\Lambda_\mu$-cb maps for $1\leq i,j\leq n$. Then the map $\phi:E\to M_n(F)$ defined as $\phi(e)=(\phi_{ij}(e))$ is $\Lambda_\mu$-cb.
\begin{proof}
Consider $\phi_X^\mu:E\bigotimes_\mu X\to M_n(F)\bigotimes_\mu X$ defined on rank one tensors as 
\[\phi_X^\mu(e\otimes x)=(\phi_{ij}(e))\otimes x=(\phi_{ij}(e)\otimes x)=\left(\phi_{ij,X}^\mu(e\otimes x)\right).\] 
Thus if $u=\sum_{i=1}^k e_i\otimes x_i\in E\bigotimes X$, we have
\[
\Vert \phi_X^\mu (u)\Vert=\Vert \phi_{ij, X}^\mu (u)\Vert\leq n^2\max_{i,j}\Vert\phi_{ij,X}^\mu\Vert\Vert u\Vert\leq n^2\max_{i,j}\Vert\phi_{ij}\Vert_\Lambda^\mu\Vert u\Vert,
\]
i.e. $\phi$ is $\Lambda_\mu$-cb with $\Vert \phi\Vert_\Lambda^\mu\leq n^2\max_{1\leq i,j\leq n}\Vert\phi_{ij}\Vert_{cb}^{\lambda_\mu}$.
\end{proof}
\end{lem}

Let us identify $M_n(\mathcal{CB}_\Lambda^\mu(E,F))$ with $\mathcal{CB}_\Lambda^\mu(E,M_n(F))$ using the bijection $\Phi:M_n(\mathcal{CB}_\Lambda^\mu(E,F))\to \mathcal{CB}_\Lambda^\mu(E,M_n(F))$, defined as $\Phi([\phi_{ij}])(e):=[\phi_{ij}(e)]$ for $e\in E$ and $\phi_{ij}\in \mathcal{CB}_\Lambda^\mu(E,F)$, $1\leq i,j\leq n$.
Instead of directly verifying the Ruan's axioms, we can easily see that there is a natural completely isometric embedding of $\mathcal{CB}_\Lambda^\mu(E,F)$ into a product space. 
Let $D_\mathcal{A}=\{u\in E\bigotimes^\mu X\mid \Vert u\Vert\leq 1\}$ and for $X\in\Lambda$ and $u\in D_X$ let $F_u=F$. 
Define $\Phi:\mathcal{CB}_\Lambda^\mu(E,F)\to \prod_{X\in\Lambda}\prod_{u\in D_X} (F_u\bigotimes^\mu X)$ as $\phi\mapsto (\phi_X^\mu(u))$.

\begin{prop}
The map $\Phi$ defined above is a complete isometry. 
\begin{proof}
Let $\phi=[\phi_{ij}]_{ij}\in M_n(\mathcal{CB}_\Lambda^\mu(E,F))=\mathcal{CB}_\Lambda^\mu(E,M_n(F))$. We have 
\[\Vert\phi\Vert=\sup_{X\in\Lambda}\Vert\phi_X^\mu\Vert
=\sup_{u\in D_X,X\in\Lambda}\Vert\phi_X^\mu(u)\Vert
=\sup_{u\in D_X,X\in\Lambda}\Vert\left[(\phi_{ij}\otimes I_X)(u)\right]_{ij}\Vert.
\]
Now consider $\Phi(\phi)=\left((\phi\otimes I_X)(u)\right)_{u,X}=\left([(\phi_{ij}\otimes I_X)(u)]_{ij}\right)_{u,X}$. Hence, as we have the sup norm on a product space, we get
\[
\Vert\Phi(\phi)\Vert=\sup_{u,X}\Vert[(\phi_{ij}\otimes I_X)(u)]_{ij}\Vert.
\]
This proves that $\Phi$ is a complete isometry.
\end{proof}
\end{prop}

The above proposition shows that $\mathcal{CB}_\Lambda^\mu(E,F)$ is an operator space equipped with natural matrix norms. The dual $E^*=\mathcal{CB}(E,\mathbb{C})$ of an operator space $E$ has a natural operator space structure inherited from the identification $M_n(\mathcal{CB}(E,\mathbb{C}))=\mathcal{CB}(E,M_n)$ \cite{Ef_Ru}. 
On the other hand, the matrix norms defined on $\mathcal{CB}_{\Lambda}(E,\mathbb{C})$ is using the identification, $M_n(\mathcal{CB}_\Lambda(E,\mathbb{C}))=\mathcal{CB}_\Lambda(E,M_n)$. 
We denote by $E^*_\Lambda$ the operator space $\mathcal{CB}_\Lambda(E,\mathbb{C})$. 
Any bounded linear functional on $E$ being completely bounded with $\Vert\phi\Vert_{cb}=\Vert\phi\Vert$, we have, $\Vert\phi\Vert_{cb}^\Lambda=\Vert\phi\Vert_{cb}$.
Hence $E^*$ and $E^*_\Lambda$ are not completely isometric in general, though they are always isometric. 
However, since $\Vert \psi\Vert_{cb}^\Lambda\leq \Vert\psi\Vert_{cb}$ for any $\psi:E\to M_n$, we have, $\Vert\psi\Vert_{M_n(E^*_\Lambda)}\leq\Vert\psi\Vert_{M_n(E^*)}$ for every $\psi\in M_n(E^*)$.

When we consider $\Lambda=\{M_n\mid n\in \mathbb{N}\}$, then the natural map gives the isometric identity $\mathcal{CB}_\Lambda(E,M_n)=\mathcal{CB}(E,M_n)$ for any $n\in \mathbb{N}$, which further gives rise to a complete isometry between $E^*$ and $E^*_\Lambda$. 
\begin{eg}
When we take $\Lambda=\Lambda_n$ as in Example \ref{eg_ckmn}, then $E^*$ and $E^*_\Lambda$ fail to be completely isometric in general. 
For example, taking $E=M_k$ for some $k>n$ and functionals $\psi_{ij}$ as the projections to $ji^{\text{th}}$ coordinate gives rise to the map $\psi=(\psi_{ij})\in \mathcal{CB}_\Lambda(E,M_k)$ as the transpose map for which $\Vert\psi\Vert_{cb}^\Lambda<\Vert\psi\Vert_{cb}$.
\end{eg}

\section{Bilinear maps associated with weighted cb maps}
\label{sec_bil}
A. Defant and D. Wiesner \cite{AD_DW} introduced a generalized way of defining operator space tensor products using special kind of multilinear maps and compared it with homogeneous polynomials. The operator space projective, Haagerup and Schur tensor products came out as three special cases in their construction.
Later, matrix ordering and related properties associated with it have been studied in \cite{PL_AK_VR}.
The following definition of $\lambda_\mu$-cb bilinear maps is motivated from their work.
Let $E,F$ and $G$ be operator spaces, and $\bigotimes^\mu$ be an operator space matrix cross tensor product.
Let $\lambda=(\lambda_n:M_n\times M_n\to M_{k(n)})$ be a sequence of bilinear maps where $n,k(n)\in \mathbb{N}$. 
We call a bilinear map $\phi:E\times F\to G$ as \emph{completely $\lambda_\mu$-bounded} (\emph{$\lambda_\mu$-cb} in short) if the associated bilinear maps $\phi_n:E\bigotimes_\mu M_n\times F\bigotimes_\mu M_n\to G\bigotimes_\mu M_{k(n)}$ defined on rank one tensors as $\phi_n(e\otimes a,f\otimes b)=\phi(e,f)\otimes \lambda_n(a,b)$ are uniformly bounded, i.e. $\sup_{n\in \mathbb{N}}\Vert \phi_n\Vert<\infty$. 

Let us denote by $\mathrm{CB}_\lambda^\mu(E\times F,G)$ the set of all $\lambda_\mu$-cb bilinear maps from $E\times F$ to $G$. 
It is easy to see that $\mathrm{CB}_\lambda^\mu(E\times F,G)$ forms a vector space and the function $\Vert\cdot\Vert_\lambda^\mu:\mathrm{CB}_\lambda^\mu(E\times F,G)\to \mathbb{R}$ defined as $\Vert\phi\Vert_\lambda^\mu=\sup_{n\in \mathbb{N}}\Vert\phi_n\Vert$ is a seminorm on it, and whenever $\lambda_n$ is nonzero for some $n\in \mathbb{N}$, it becomes a norm. Hence we shall always consider the case when not all $\lambda_n$ are zero bilinear maps and we assume further that $\lambda$ is uniformly bounded.
Now let us give an operator space structure to $\mathrm{CB}_\lambda^\mu(E\times F,G)$ using the (natural) identification $M_n(\mathrm{CB}_\lambda^\mu(E\times F,G))=\mathrm{CB}_\lambda^\mu(E\times F,M_n(G))$.

\begin{thm} If $E,F,G$ are operator spaces and $\phi^{ij}\in \mathrm{CB}_\lambda^\mu(E\times F,G)$ for $1\leq i,j\leq m$, then the bilinear map $\phi:E\times F\to M_m(G)$ defined as $\phi(e,f)=[\phi^{ij}(e,f)]$ is completely $\lambda_\mu$-bounded with $\Vert\phi\Vert_\lambda^\mu\leq n^2\max_{i,j}\Vert\phi^{ij}\Vert_\lambda^\mu$. Moreover, $\mathrm{CB}_\lambda^\mu(E\times F,G)$ is an operator space with matrix norms induced from the identification $M_n(\mathrm{CB}_\lambda^\mu(E\times F,G))=\mathrm{CB}_\lambda^\mu(E\times F,M_n(G))$.
\begin{proof}
Consider $\phi_n:E\bigotimes_\mu M_n\times F\bigotimes_\mu M_n\to M_m(G)\bigotimes_\mu M_{k(n)}$. 
Let $u_1=\sum_{r=1}^s e_r\otimes a_r\in E\bigotimes M_n$ and $u_2=\sum_{r=1}^s f_r\otimes b_r\in F\bigotimes M_n$. Then 
$\phi_n(u_1,u_2)=[\sum_{r=1}^s\phi^{ij}(e_r,f_r)\otimes \lambda_n(a_r,b_r)]$. Thus 
\begin{equation*}
\begin{aligned}
\Vert\phi_n(u_1,u_2)\Vert &= \Vert [\phi_n^{ij}(u_1,u_2)]\Vert\\
&\leq m^2 \max_{i,j}\Vert\phi_n^{ij}(u_1,u_2)\Vert\\
&\leq m^2 \max_{i,j}\Vert\phi_n^{ij}\Vert\Vert u_1\Vert\Vert u_2\Vert\\
&\leq m^2 \max_{i,j}\Vert\phi^{ij}\Vert_\lambda^\mu\Vert u_1\Vert\Vert u_2\Vert,
\end{aligned}
\end{equation*}
i.e. $\phi= [\phi^{ij}]\in \mathrm{CB}_\lambda^\mu(E\times F,M_m(G))$ and $\Vert\phi\Vert_\lambda^\mu\leq m^2\max_{i,j}\Vert\phi^{ij}\Vert_\lambda^\mu$.

For the other part, we shall verify the Ruan's axioms in order to prove that it is an operator space.

Let $\phi_1\in M_n(\mathrm{CB}_\lambda^\mu(E\times F,G))=\mathrm{CB}_\lambda^\mu(E\times F,M_n(G))$ and $\phi_2\in M_m(\mathrm{CB}_\lambda^\mu(E\times F,G))=\mathrm{CB}_\lambda^\mu(E\times F,M_m(G))$. 
Consider $\phi_1\oplus\phi_2\in M_{n+m}(\mathrm{CB}_\lambda^\mu(E\times F,G))=\mathrm{CB}_\lambda^\mu(E\times F,M_{n+m}(G))$. 
Let $u_1=\sum_{r=1}^s e_r\otimes a_r\in E\bigotimes M_l$ and $u_2=\sum_{r=1}^s f_r\otimes b_r\in F\bigotimes M_l$. 
Then
\begin{equation*}
\begin{aligned}
\Vert(\phi^1\oplus\phi^2)_l(u_1,u_2)\Vert &= \left\Vert \begin{bmatrix}
\phi_l^1(u_1,u_2) & 0\\
0 & \phi_l^2(u_1,u_2) 
\end{bmatrix}
\right\Vert\\
&\leq \max_{i=1,2}\Vert\phi_l^i(u_1,u_2)\Vert\\
&\leq \max_{i=1,2}\Vert\phi_l^i\Vert\Vert u_1\Vert\Vert u_2\Vert\\
&\leq \max_{i=1,2}\Vert\phi^i\Vert_\lambda^\mu\Vert u_1\Vert\Vert u_2\Vert.
\end{aligned}
\end{equation*}
Thus $\Vert\phi^1\oplus\phi^2\Vert_\lambda^\mu\leq \max_{i=1,2}\Vert\phi^i\Vert_\lambda^\mu$. 

Let $\phi=[\phi^{ij}]\in M_n(\mathrm{CB}_\lambda^\mu(E\times F,G))$ and $\alpha,\beta$ be scalar rectangular matrices such that $\alpha\phi\beta\in M_m(\mathrm{CB}_\lambda^\mu(E\times F,G))=\mathrm{CB}_\lambda^\mu(E\times F,M_m(G))$. 
Then
\begin{equation*}
\begin{aligned}
\Vert(\alpha\phi\beta)_l(u_1,u_2)\Vert &= \Vert (\alpha[\phi^{ij}]\beta)_l(u_1,u_2)\Vert\\
&=\Vert \alpha[\phi_l^{ij}(u_1,u_2)]\beta\Vert\\
&\leq\Vert\alpha\Vert\Vert[\phi_l^{ij}(u_1,u_2)]\Vert\Vert\beta\Vert\\
&\leq\Vert\alpha\Vert\Vert\phi_l\Vert\Vert u_1\Vert\Vert u_2\Vert\Vert\beta\Vert\\
&\leq\Vert\alpha\Vert\Vert\phi\Vert_\lambda^\mu\Vert\Vert\beta\Vert\Vert u_1\Vert\Vert u_2\Vert.
\end{aligned}
\end{equation*}
Thus $\Vert\alpha\phi\beta\Vert_\lambda^\mu\leq \Vert\alpha\Vert\Vert\phi\Vert_\lambda^\mu\Vert\Vert\beta\Vert$.

As $\Vert\cdot\Vert_\lambda^\mu$ is a norm on $M_m(\mathrm{CB}_\lambda^\mu(E\times F,G))$ for all $m\in \mathbb{N}$, it follows that the above two inequalities are sufficient to conclude that $\mathrm{CB}_\lambda^\mu(E\times F,G)$ is an operator space.
\end{proof}
\end{thm}

\begin{rem}
\label{bil_intertwine}
It is very important to observe that intertwining the roles of $E$ and $F$ in $\mathrm{CB}_\lambda^\mu(E\times F,G)$ matters, i.e. in general, $\mathrm{CB}_\lambda^\mu(E\times F,G)$ and $\mathrm{CB}_\lambda^\mu(F\times E,G)$ are not completely isometric. For example, consider $\lambda=(\lambda_n:M_n\times M_n\to M_n)$ defined as $\lambda_n(a,b)=ab$. Then the non-commutativity of the matrix multiplication plays an important role in making the difference.  
\end{rem}

We shall now consider three very special cases as follows:
\begin{enumerate}[\text{Case} 1:]
\item $\lambda=(\lambda_n:M_n\times M_n\to M_n)$ defined as $\lambda_n(a,b)=ab$, the usual matrix multiplication representing composition of linear maps.
\item $\lambda=(\lambda_n:M_n\times M_n\to M_{n^2})$ defined as $\lambda_n(a,b)=a\otimes b$, where $a\otimes b$ denotes the Kronecker multiplication, i.e. for $a=[a_{ij}],b=[b_{ij}]\in M_n$ the pair $(a,b)$ is mapped to $[a_{ij}b_{kl}]$.
\item $\lambda=(\lambda_n:M_n\times M_n\to M_n)$ defined as $\lambda_n(a,b)=a\odot b$, where $a\odot b$ denotes the Schur multiplication of two matrices, i.e. for $a=[a_{ij}],b=[b_{ij}]\in M_n$ the pair $(a,b)$ is mapped to $[a_{ij}b_{ij}]$.
\end{enumerate}

Clearly, when we take $\mu=min$, then the above three cases give rise to the corresponding spaces of bilinear maps $\mathrm{CB}_\lambda^\mu(E\times F,G)$ as the spaces of matricially completely bounded bilinear maps, jointly completely bounded bilinear maps, and completely Schur bounded bilinear maps respectively. We shall have a closer look at each of them.

When we take $G=\mathbb{C}$, we call $\psi\in \mathrm{CB}_\lambda^\mu(E\times F,\mathbb{C})$ as a \emph{$\lambda_\mu$-cb bilinear form}. 
In the same way in which we associate jointly completely bounded (respectively completely bounded) bilinear forms with the operator space projective (respectively Haagerup) tensor products \cite{Bl_Pn}, we wish to associate a tensor norm $\Vert\cdot\Vert_\lambda^\mu$ on the algebraic tensor product $E\bigotimes F$ using the completely $\lambda_\mu$-bounded bilinear forms on $E\times F$ such that the operator space dual of the completed tensor product $E\bigotimes^{\lambda_\mu} F$ is $\mathrm{CB}_\lambda^\mu(E\times F,\mathbb{C})$ completely isometrically, 
i.e. we shall define $\Vert\cdot\Vert_\lambda^\mu$ on $E\bigotimes F$ as follows: for $u\in M_n(E\bigotimes F)$, $\Vert u\Vert_\lambda^\mu=\sup\{\Vert\llangle u,\psi\rrangle\Vert\mid \psi\in M_m(\mathrm{CB}_\lambda^\mu(E\times F,\mathbb{C})),\Vert \psi\Vert_\lambda^\mu\leq 1\}$ where $\llangle u,\psi\rrangle$ denotes the matrix pairing with the abuse of notation of identifying linear maps on $E\bigotimes F$ with bilinear maps on $E\times F$. 
If the quantity $\Vert u\Vert_\lambda^\mu$ is finite for every $u\in M_n(E\bigotimes F)$, then clearly it defines a norm on $M_n(E\bigotimes F)$ for every $n$, and the matrix norms satisfy the Ruan's conditions for an operator space norm.
The operator space structure obtained by doing so is called the \emph{dual operator space structure}. We denote by $E\bigotimes^{\lambda_\mu} F$ the completion with respect this norm, which we call the $\lambda_\mu$ tensor product of $E$ and $F$.

In the three cases mentioned above, we need to verify whether $\sup\{\Vert\llangle u,\psi\rrangle\Vert\mid \psi\in M_m(\mathrm{CB}_\lambda^\mu(E\times F,\mathbb{C})),\Vert \psi\Vert_\lambda^\mu\leq 1\}$ is finite for each $u\in M_n(E\bigotimes F)$. 
Let us denote by $\bigotimes_{\lambda_\mu}$ the natural map $E\bigotimes_\mu M_n\times F\bigotimes_\mu M_n\to E\bigotimes F\bigotimes M_{k(n)}$ given by $\bigotimes_{\lambda_\mu}(e\otimes a,f\otimes b)=a\otimes f\otimes \lambda(a,b)$.

In case 1, let $u\in M_n(E\bigotimes F)$. Then we can write $u$ in the form $u=\bigotimes_{\lambda_\mu}(e,f)$ for some $e=[e_{ij}]\in M_m^\mu(E)$ and $f=[f_{kl}]\in M_m^\mu(F)$. 
Let $\phi=[\phi_{pq}]\in M_p(\mathrm{CB}_\lambda^\mu(E\times F,\mathbb{C}))$.
Consider 
\begin{equation*}
\begin{aligned}
\left\llangle \bar{\phi},\textstyle\bigotimes_{\lambda_\mu}(e,f)\right\rrangle &= \left\llangle \left[\bar{\phi}_{pq}\right]_{p,q},\left[\sum_{j=1}^m e_{ij}\otimes f_{jl}\right]_{i,l} \right\rrangle\\
&= \left[\bar{\phi}_{pq}\left(\sum_{j=1}^m e_{ij\otimes f_{jl}}\right)\right]_{p,q,i,l}\\
&= \left[(\phi_{ij})_n(e,f)\right]\\
&= U\phi_n(e,f)V,
\end{aligned}
\end{equation*}
for some invertible scalar matrices $U,V$ such that $\Vert U\Vert=1=\Vert V\Vert$. 
Hence \[
\Vert\llangle\phi,\textstyle\bigotimes_{\lambda_\mu}(e,f)\rrangle\Vert\leq \Vert U\phi_n(e,f)V\Vert\leq \Vert\phi\Vert_\lambda^\mu\Vert e\Vert_{M_n^\mu{E}}\Vert f\Vert_{M_n^\mu(F)}.
\]

A similar calculation shows the finiteness of $\sup\{\Vert\llangle u,\psi\rrangle\Vert\mid \psi\in M_m(\mathrm{CB}_\lambda^\mu(E\times F,\mathbb{C})),\Vert \psi\Vert_\lambda^\mu\leq 1\}$ for any $u\in M_n(E\bigotimes F)$ in the other two cases also.

Deducing even basic properties of the tensor product such as associativity, functorial property are extremely challenging, though we can obtain commutativity of the tensor product by putting reasonable restrictions on $\lambda$.
We say $\lambda$ is \emph{symmetric} if there exists a sequence of unitary matrices $(u_k)$ where $u_k\in M_k$ such that $\lambda_n(b,a)=u_{k(n)}^{-1}\lambda_n(a,b)u_{k(n)}$ for all $a,b\in M_n$, $n\in \mathbb{N}$. In Case 2 and Case 3 above, the given $\lambda$ is symmetric. 
Observe that $\Vert\lambda_n(b,a)\Vert=\Vert u_{k(n)}^{-1}\lambda_n(a,b)u_{k(n)}\Vert$ and hence it clearly follows that whenever $\lambda$ is symmetric, $\mathrm{CB}_\lambda^\mu(E\times F,G)=\mathrm{CB}_\lambda^\mu(F\times E,G)$ completely isometrically via the natural map $\Phi$ defined as $\Phi(\phi)=\bar{\phi}$ where $\bar{\phi}(f,e)=\phi(e,f)$.

\begin{thm}
If $\lambda$ is symmetric, then $\bigotimes^{\lambda_\mu}$ is commutative, i.e. for any two operator spaces $E,F$ we have completely isometrically $E\bigotimes^{\lambda_\mu} F=F\bigotimes^{\lambda_\mu} E$.
\begin{proof}
As $\lambda$ is symmetric, we have completely isometrically $\mathrm{CB}_\lambda^\mu(E\times F,\mathbb{C})\stackrel{\Phi}{=} \mathrm{CB}_\lambda^\mu(F\times E,\mathbb{C})$. 
But $\mathrm{CB}_\lambda^\mu(E\times F,\mathbb{C})$ is nothing but the dual of $E\bigotimes^{\lambda_\mu} F$, and $\mathrm{CB}_\lambda^\mu(F\times E,\mathbb{C})$ is that of $F\bigotimes^{\lambda_\mu} E$. 
For $\phi\in\mathrm{CB}_\lambda^\mu(E\times F,\mathbb{C})$, we have 
\[
\langle\phi, e\otimes f\rangle=\phi(e,f)=\Phi(\phi)(f,e)=\langle\Phi(\phi),f\otimes e\rangle.
\]
Hence the linear isomorphism $\Psi:E\bigotimes F\to F\bigotimes E$ defined on elementary tensors as $\Psi(e\otimes f)=f\otimes e$ extends to a complete isometry between $E\bigotimes^{\lambda_\mu} F$ and $F\bigotimes^{\lambda_\mu} E$, because if $(u_k)$ is a Cauchy sequence in $M_n(E\bigotimes_\lambda F)$, then $(\Psi^{(n)}(u_k))$ is Cauchy in $M_n(F\bigotimes_\lambda E)$ with $\Vert u_k\Vert=\Vert\Psi^{(n)}(u_k)\Vert$. 
\end{proof}
\end{thm}

As we have mentioned in Section \ref{subsec_quant}, we could not find a reasonable sufficient condition for a tensor product to be quantizing. It would be interesting if one can completely characterize the quantizing tensor products. 
An analogue of $\lambda_\mu$-cb maps when $\mu$ being non-quantizable and additionally considering a random class $\Lambda$ of Banach spaces/operator spaces as discussed in Section \ref{sec_lambda-cb} will be a generalization the ideas that we discussed, though it can be more complicated.

\section*{Acknowledgements}
Research of the first author is supported by the Council of Scientific \& Industrial Research (CSIR), Govt. of India (Ref. No.: 09/045(1403)/2015-EMR-I).

\bibliographystyle{amsplain}

\end{document}